\documentclass{amsart}
\usepackage{amsmath,amsfonts,amssymb}
\usepackage{graphicx}

\newtheorem{theorem}{Theorem}[section]

\newtheorem{cor}[theorem]{Corollary}
\newtheorem{THEO}{Theorem}

\theoremstyle{definition}

\theoremstyle{remark}

\numberwithin{equation}{section}

\begin{document}

\title[Sharp Hardy's Inequalities in Hilbert Spaces]{Sharp Hardy's Inequalities in Hilbert Spaces}
\thanks{Research supported by the Brazilian Science Foundations FAPESP under Grants 2016/09906-0 and  2021/13340-0 and CNPq under Grant 309955/2021-1, and the 
Bulgarian National Research Fund through Contract KP-06-N62/4}

\author[D. K. Dimitrov]{Dimitar K. Dimitrov}
\address{Departamento de Matem\'atica, IBILCE, Universidade Estadual Paulista, 15054-000
S\~{a}o Jos\'{e} do Rio Preto, SP, Brazil }
\email{d\_k\_dimitrov@yahoo.com}

\author[Ivan Gadjev]{Ivan Gadjev}
 \address{Department of Mathematics and Informatics,
University of Sofia,
5 James Bourchier Blvd., 
1164 Sofia, Bulgaria}
\email{gadjev@fmi.uni-sofia.bg}

\author[M. E. H. Ismail]{Mourad E. H. Ismail}
\address{Department of Mathematics, University of Central Florida, Orlando, Florida 32816, USA}
%\curraddr{}
\email{ismail@math.ucf.edu}

      

\begin{abstract}
We study the behavior of the smallest possible constants $d(a,b)$ and $d_n$   in Hardy's inequalities
$$
\int_a^b\left(\frac{1}{x}\int_a^xf(t)dt\right)^2\,dx\leq
d(a,b)\,\int_a^b [f(x)]^2 dx
$$
and
$$
\sum_{k=1}^{n}\Big(\frac{1}{k}\sum_{j=1}^{k}a_j\Big)^2\leq
d_n\,\sum_{k=1}^{n}a_k^2.
$$
The exact constant $d(a,b)$ and the precise rate of convergence of $d_n$ are established and the extremal function 
and  the ``almost extremal'' sequence are found.
\end{abstract}

\maketitle


%==================================================================
\section{Introduction and statement of the results}
%==================================================================

In the series of papers  \cite {H1919, H1920, H1925} Hardy proved the following two inequalities. Let $p>1$.
If $f(x)\geq0$ and $f^p$ is integrable over $(0,\infty)$, then the inequality 
\begin{equation}\label{eq11}
\int_0^\infty\left(\frac{1}{x}\int_0^xf(t)dt\right)^p\,dx\leq
\left(\frac{p}{p-1}\right)^{p}\,\int_0^\infty [f(x)]^p\, dx.
\end{equation}
holds.
This is the original Hardy's integral inequality. 
The discrete version of Hardy's inequality reads as
\begin{equation}\label{eq12}
\sum_{k=1}^{\infty}\Big(\frac{1}{k}\sum_{j=1}^{k}a_j\Big)^p\leq
\left(\frac{p}{p-1}\right)^{p}\,\sum_{k=1}^{\infty}{a_k^p},\ \ \  \ a_k \geq 0,\quad k\in \mathbb{N}.
\end{equation}
Initially Hardy proved (\ref{eq12}) with the constant $(p^2/(p-1))^p$.
Later Landau, in the letter \cite{Lan2},  which was officially published in \cite{Lan1}, 
established  the exact constant $(p/(p-1))^{p}$ in the sense that 
there is no smaller one for which \eqref{eq12}  holds for every sequence of nonnegative numbers $a_k$. Moreover, 
Landau observed that equality in \eqref{eq12} occurs only for the trivial sequence, that is,
 when $a_k=0$ for every $k \in \mathbb{N}$. Similarly, equality in \eqref{eq11} occurs if and only
if $f(x)\equiv0$ almost everywhere.

Inequality \eqref{eq11} has been extended to what is nowadays called the {\em general Hardy integral inequality}: 
\begin{equation}\label{eq15}
\int_a^b\left(\frac{1}{x}\int_a^xf(t)dt\right)^p\,dx\leq
d_p(a,b)\,\int_a^b [f(x)]^p dx, \qquad f(x)\geq0.
\end{equation}

For a fixed $n \in \mathbb{N}$, it is natural to study the inequality 
\begin{equation}\label{eq16}
\sum_{k=1}^{n}\Big(\frac{1}{k}\sum_{j=1}^{k}a_j\Big)^p\leq
d_{n,p} \,\sum_{k=1}^{n}{a_k^p}, \qquad a_k\geq0,\,\,\,k=1,2,...,n.
\end{equation}
and ask for the smallest possible $d_{n,p}$, for which it holds.

When $p$ is a positive even integer the assumption for nonnegativity of $f(x)$ and $\,\{a_k\}\,$  
can be dropped. In particular, when $p=2$ the inequalities \eqref{eq15} and \eqref{eq16}  become
\begin{equation}\label{eq13}
\int_a^b\left(\frac{1}{x}\int_a^xf(t)dt\right)^2\,dx\leq
d(a,b)\,\int_a^b [f(x)]^2\, dx
\end{equation}
and 
\begin{equation}
\label{eq14}
\sum_{k=1}^{n}\Big(\frac{1}{k}\sum_{j=1}^{k}a_j\Big)^2\leq
d_n\,\sum_{k=1}^{n}{a_k^2}.
\end{equation}
 
There are many papers investigating different generalizations and applications of Hardy's inequality;  
see for instance \cite{KMP2007} and the bibliography in the book \cite{KP2003}.

In the present paper we establish the sharp inequality (\ref{eq13}) in the sense that we determine the exact constant $d(a,b)$ in \eqref{eq13} as well as an extremal function $f$, for which equality is attained. Our main result reads as follows: 

\begin{theorem}\label{th1}
Let  $a$ and $b$ be any fixed numbers with $0<a<b<\infty$. Then the inequality
\begin{equation}\label{MR1}
\int_a^b\left(\frac{1}{x}\int_a^xf(t)dt\right)^2 dx\leq
\frac{4}{1+4\alpha^2}\,\int_a^b [f(x)]^2 \, dx,
\end{equation}
where $\alpha$ is the only solution of the equation
$$
\tan\left(\alpha\log\frac{b}{a} \right)+2\alpha=0\ \ \mathrm{in\ the\ interval}\ \ \left(\frac{\pi}{2\log\frac{b}{a}},\frac{\pi}{\log\frac{b}{a}} \right),
$$
holds for every $f \in L^2[a,b]$.
Moreover, equality in \eqref{MR1} is attained for 
\[
f_{a,b}(x)=x^{-1/2}\left(2\alpha\cos(\alpha\log x)+\sin(\alpha\log x) \right).
\]
\end{theorem}

\begin{cor}\label{cor1}
 When either of the limits relations $a\rightarrow0$,  $b\rightarrow\infty$, or both hold, i.e. $\log(b/a)\rightarrow\infty$, then
\[
d(a,b) \sim 
4-\frac{c}{[\log\frac{b}{a}]^2}.
\]
In other words, there exist absolute constants $c_1>0$ and $c_2>0$, such that
\begin{equation*}
4-\frac{c_1}{[\log\frac{b}{a}]^2}\le d(a,b) \le 
4-\frac{c_2}{[\log\frac{b}{a}]^2}.
\end{equation*}
%and it is obvious that when $b\rightarrow\infty$  the function $f_{a,b}$ converges uniformly to the function $g(x)\equiv 0$.
\end{cor}
Since for the function $f_{a,b}$ defined above, obeys
\[
0<f_{a,b}(x)<x^{-1/2}\alpha\left(2+|\ln x| \right)
<\frac{\pi}{\ln\frac{b}{a}}x^{-1/2}\left(2+|\ln x| \right),  
\] 
it is obvious that $f_{a,b}$ converges uniformly to $g(x)\equiv 0$ when $b\rightarrow\infty$.

An important well-known fact is that (\ref{eq11}) is the prototype of the Hardy-Littlewood inequality for the maximal function. 
Therefore, it is natural to consider the inequality
\begin{equation*}
\int_a^b\left(\frac{1}{x-a}\int_a^xf(t)dt\right)^2\,dx\leq
d(a,b)\,\int_a^b [f(x)]^2 dx, \ \ f \in L^2[a,b]
\end{equation*}
which is equivalent, by change of variables, to
\begin{equation*}
\int_0^b\left(\frac{1}{x}\int_0^xf(t)dt\right)^2\,dx\leq
d(b)\,\int_0^b [f(x)]^2 dx, \ \ f \in L^2[0,b].
\end{equation*}
Then, Corollary \ref{cor1} implies that $d(b)=4$ and the only function $f$ for which the equality is attained is the one which vanishes almost everywhere.

Inequality (\ref{eq14}) has been studied by many authors. Based on ideas of Widom \cite{Widom} and his previous joint work with Widom himself \cite{WW}, and with N. de Bruijn \cite{deBW}, on Hilbert's inequality, Herbert Wilf 
\cite {Wilf} established the following asymptotic expression for $d_n$:
\begin{equation}
\label{WA}
d_n = 4 - \frac{16 \pi^2}{\log^2 n} + O \left( \frac{\log \log n}{\log^3 n} \right),\ \ \ n \to \infty.
\end{equation}
The ideas developed in \cite{DIGR} and an important observation of the third-named author of the present note, yielded an explicit representation of $d_n$ in terms of the smallest zero of a continuous dual Hahn polynomial of degree $n$, for a specific choice of the parameters, in terms of which these polynomials are defined; see Theorem A below.  More recently F. Stampach \cite{Stampach} studied in details the asymptotics of $d_n$, again employing its relation to these zeros and proved that    
\begin{equation}
\label{SA}
d_n = 4 - \frac{16 \pi^2}{\log^2 n} + \frac{32 \pi^2 (\gamma+ 6 \log 2)}{\log^3 n} + O \left( \frac{\log \log n}{\log^4 n} \right),\ \ \ n \to \infty.
\end{equation}
Another formula involving few more terms and an algorithm for calculating the further ones in the asymptotic expansion of $d_n$ in terms of negative powers of $\log n$ was suggested in \cite{Stampach}. 

Here we combine an idea similar to the one we use for the proof of Theorem \ref{th1} with the result in \cite[Theorem 1.1]{DIGR} concerning the relation between $d_n$ and the zeros of the continuous dual Hahn polynomial, and a recent one, due to W-G. Long, D. Dai, Y-T. Li, X-S. Wang \cite{Long} about the asymptotic behaviour of those zeros. We establish sharp lower and upper bounds for $d_n$ and obtain its full asyptotic expansion, thus extending the results of Wilf (\ref{WA}) and  Stampach (\ref{SA}) to the largest possible generality.

\begin{theorem}\label{th2}
Let
\[
a_k=\int_k^{k+1}h(x)dx,
\] 
where 
\begin{equation}\label{eq20}
h(x)=x^{-1/2}\left(2\alpha\cos(\alpha\log x)+\sin(\alpha\log x) \right),\quad 1\le x\le n+1,
\end{equation}
and $\alpha$ is the only solution of the equation
\[
\tan(\alpha\log (n+1))+2\alpha=0\ \ \mathrm{in\ the\ interval}\ \ \left(\frac{\pi}{2\log(n+1)},\frac{\pi}{\log(n+1)} \right).
\]
 Then 
\begin{equation}\label{eq18}
\sum_{k=1}^{n}\Big(\frac{1}{k}\sum_{j=1}^{k}a_j\Big)^2\geq
\frac{4}{1+4\alpha^2}\,\sum_{k=1}^{n}{a_k^2}.
\end{equation}
 \end{theorem}

Since
\[
\frac{4}{1+4\alpha^2}\geq 4-16\alpha^2>4-\frac{16\pi^2}{\log^2(n+1)}
\]
we obtain 
\[
d_n\geq 4-\frac{16\pi^2}{\log^2(n+1)}.
\]

Combining the above mentioned results with the latter we obtain the exact rate of convergence of $\{ d_n \}$ as well as very sharp estimates for $d_n$ for every fixed $n$:
\begin{theorem}\label{Th1.4}
The inequalities 
\begin{equation}
\label{estdn}
4-\frac{16\pi^2}{[\log(n+1)]^2} \le d_n \le 
4-\frac{32}{[\log n + 4]^2}
\end{equation}
hold for every natural $n\geq 3$.

Moreover, 
\begin{equation}
\label{asdn}
d_n \sim 
4-\frac{16\, \pi^2}{4\, \pi^2 + (\gamma + 6 \log 2 +\log n)^2}\ \ \mathrm{as}\ \ n\rightarrow \infty,
\end{equation}
where $\gamma$ is the Euler constant, and the following full asyptotic expansion of $d_n$ in terms of the negative powers of $\log n$ holds: For every fixed $m\in \mathbb{N}$, $m \geq 2$,
\begin{equation}
\label{asexp}
d_n \sim 
4- \sum_{k=2}^m c_k\  \frac{1}{[\log n]^k}\, +\, O((\log n)^{-m-1}),\ \ \ n \to \infty, 
\end{equation}
with 
\begin{equation}
\label{ck}
c_k = 16 \pi^2  \left(4 \pi^2 + (\gamma+6\log2)^2\right)^{(k-2)/2}\ U_{k-2}\left( -\frac{\gamma+6\log2}{\sqrt{4 \pi^2 + (\gamma+6\log2)^2}} \right),
\end{equation}
where $U_{k-2}$ denotes the Chebyshev polynomial of the second kind of degree $k-2$.
\end{theorem}

Needless to say, the first few coefficients $c_k$, $k=2,3,4,5$, in (\ref{ck}) coincide with those obtained by  Wilf and Stampach. 

The lower bound in (\ref{estdn}), obtained as a consequence of Theorem \ref{th2}, is amazingly close to the asymptotic value for $d_n$ in (\ref{asdn}), which is derived via \cite[Theorem 1.1]{DIGR} and  \cite{Long}.  Indeed, it is not difficult to verify that for their difference one has 
$$
\left( 4-\frac{16\, \pi^2}{4\, \pi^2 + (\gamma + \log(64 n))^2} \right) - \left( 4-\frac{16\pi^2}{\log^2(n+1)} \right) \sim \frac{16\, \pi^2(2\gamma+ \log 128)}{\log^3 n},\ \  n\rightarrow \infty.
$$
The function $h(x)$, defined in \eqref{eq20},  obeys
\[
0<h(x)<x^{-1/2}\left(2\alpha+\alpha\log x \right)
<\frac{\pi}{\log(n+1)}x^{-1/2}\left(2+\log x \right).  
\] 
Then it is obvious, as Landau pointed in his letter to Hardy, that if we let $n\rightarrow\infty$ %in \eqref{eq18} 
the almost extremal sequence $a_k,\, k=1,2,...$ defined in the Theorem \ref{th2} goes to the zero sequence, i.e.
to the sequence $a_k=0$ for all $k$.

%==================================================================
\section{Proof of Theorem \ref{th1}}
%==================================================================
 By simple changes of variables $t=au,\,x=av$ in the left-hand side and $x=au$ in the right-hand side we write 
the inequality \eqref{MR1} in the following way
\begin{equation*}
\int_1^{b/a}\left(\frac{1}{v}\int_1^vf(au)du\right)^2 dv\leq
\frac{4}{1+4\alpha^2}\,\int_1^{b/a} [f(au)]^2\, du.
\end{equation*}
%Denoting $f(ax)=g(x)$
It is obvious (by changing the notations), that it suffices to prove Theorem \ref{th1} for the interval $(1,b)$.

Cauchy's inequality yields 
\[
\left(\int_1^x f(t)dt \right)^2   \le\left(\int_1^x [g(t)]^2\, dt  \right) \left(\int_1^x \frac{[f(t)]^2}{[g(t)]^2}\, dt  \right)
\]
for every pair of  functions $f, g \in L^2[1,b]$, such that $g(x)\neq 0$ for every $x\in(1,b)$.
By multiplying both sides of the latter inequality by $x^{-2}$ and integrating from 1 to $b$ we obtain
\begin{align*}
\int_1^b\left(\frac{1}{x}\int_1^x f(t)\, dt \right)^2dx   \le 
\int_1^b \left(\frac{1}{x^2}\int_1^x [g(t)]^2\, dt  \right) \left(\int_1^x \frac{[f(t)]^2}{[g(t)]^2}\, dt  \right)
\end{align*}
By changing the order of integration in the right-hand side,  
we obtain
\begin{align*}
\int_1^b\left(\frac{1}{x}\int_1^x f(t)\, dt \right)^2dx   \le 
\int_1^b\left(\frac{1}{[g(t)]^2}\int_t^b\left(\frac{1}{x^2}\int_1^x [g(u)]^2\, du \right)dx  \right)[f(t)]^2\, dt.
\end{align*}
Let us denote for brevity %$M(t,f)=f^{-2}(t)M^*(t,f)$ where
\[
M(g,t)=\frac{1}{g^2(t)}\int_t^b\left(\frac{1}{x^2}\int_1^x [g(u)]^2\, du \right)\, dx.
\]
Then the 
inequality 
\begin{align*}
\int_1^b\left(\frac{1}{x}\int_1^x f(t)dt \right)^2dx   \le 
\max_{1< t< b}M(g,t)  \int_1^b [f(t)]^2\, dt
\end{align*}
holds for every two functions $f, g \in L^2[1,b]$, with $g(x)\neq 0$ for $1<x<b$, and consequently 
\[
d(1,b)\le \max_{1< t< b}M(g,t)
\]
for every non-vanishing function $g\in L^2(1,b)$.

Now we minimize 
\[
\max_{1< t< b}M(g,t)
\]
 over all such functions  $g$, that is, determine
\[
\min_{g(x)\neq 0}\,\max_{1< t< b}\frac{1}{g^2(t)}\int_t^b\left(\frac{1}{x^2}\int_1^x [g(u)]^2\, du \right) dx.
\]

Let us consider the function $h(x)=x^{-1/2}\left(2\alpha\cos(\alpha\log x)+\sin(\alpha\log x) \right)$, $x\in[1,b]$
where $\alpha$ is the only solution of the equation
\[
\tan\left(\alpha\log b \right)+2\alpha=0
\]
in the interval $\left(\pi/(2\log b),\pi/\log b \right)$.
We have 
\[
h'(x)=-\frac{1+4\alpha^2}{2x^{3/2}}\sin(\alpha\ln x).
\]
Obviously, $h'(x)< 0$ for $1< x< b$ and consequently $h(x)$ is decreasing. 
Since $h(b)=0$, it follows that $h(x)>0$ for $1\le x< b$. Then for the function $g(x)=\sqrt{h(x)}$ we have
$$
\int_1^x [g(u)]^2\, du = \int_1^x u^{-1/2}\left(2\alpha\cos(\alpha\log u)+\sin(\alpha\log u)\right)du=
2\sqrt x\sin(\alpha\log x),
$$
and
$$
\int_t^b\left(\frac{1}{x^2}\int_1^x [g(u)]^2 du \right)dx=
2\int_t^b\frac{\sin(\alpha\log x)}{x^{3/2}}dx=-\frac{4}{1+4\alpha^2}\int_t^bh'(x)dx
=\frac{4\, [g(t)]^2}{1+4\alpha^2}.
$$
Therefore,
\[
M(g,t)=\frac{4}{1+4\alpha^2} \quad\mbox{and}\quad d(1,b)\le \frac{4}{1+4\alpha^2}.
\]

On the other hand, by changing the order of integration, we can rewrite the left-hand side of \eqref{MR1}, with $a=1$, in the 
following way
\begin{align*}
\int_1^b\left(\frac{1}{x}\int_1^xf(t)\, dt\right)^2 dx=
\int_1^b\left(\frac{1}{f(t)}\int_t^b\left(\frac{1}{x^2}\int_1^x f(u)\, du \right)dx  \right)[f(t)]^2\, dt.
\end{align*}
Consequently, for the function $f_{a,b}(x)=h(x)$ we have
\begin{align*}
\int_1^b\left(\frac{1}{x}\int_1^xf_{a,b}(t)dt\right)^2 \,dx=
\frac{4}{1+4\alpha^2}\,\int_1^b [f_{a,b}(x)]^2\, dx.
\end{align*}
The proof of Theorem \ref{th1} is complete.

%\end{proof} 

%==================================================================
\section{Proof of Theorem \ref{th2}}
%==================================================================

%Motivated by the question about the existence of ``almost extremizing'' , 
%in the present paper we examine the best constant $d_n$  in the

%\begin{proof}

By changing the order of summation in the left hand side of \eqref{eq18} we obtain
\begin{equation*}
\sum_{k=1}^{n}\left(\frac{1}{k}\sum_{j=1}^{k}a_j\right)^2=
\sum_{i=1}^{n} \left[\frac{1}{a_i}  \left(\sum_{k=i}^{n}\frac{1}{k^2}\sum_{j=1}^{k}a_j\right)\right]a_i^2
=\sum_{i=1}^{n}M_i\, a_i^2
\end{equation*}
where
\[
M_i=\frac{1}{a_i} M_i^* \quad\mbox{and}\quad M_i^*=\sum_{k=i}^{n}\left(\frac{1}{k^2}\sum_{j=1}^{k}a_j\right).
\]
Then
\begin{equation*}
\sum_{k=1}^{n}\left(\frac{1}{k}\sum_{j=1}^{k}a_j\right)^2\geq \min_{1\le i\le n}M_i\sum_{k=1}^{n}a_i^2\, ,
\end{equation*}
and consequently
\begin{equation}\label{eq3}
d_n\geq \min_{1\le i\le n}M_i.
\end{equation}
For the function $h(x)$ defined in \eqref{eq20}  we have 
\[
h'(x)=-\frac{1+4\alpha^2}{2x^{3/2}}\sin(\alpha\log x),
\]
Obviously $h'(x)<0$ for $1<x\le (n+1)$, that is, $h(x)$ is decreasing. Since $h(n+1)=0$ it follows that $h(x)>0$.

We shall prove that the sequence 
\[
a_k=\int_k^{k+1}h(x)\, dx
\] 
is the  ``almost extremal'' sequence for Hardy's inequality \eqref{eq14}, i.e. that inequality \eqref{eq18} holds.
Since the function $h(x)$ is continuous there exists a point $\eta_i\in [i,i+1]$ such that $a_i=h(\eta_i)$.

We have
\[
\sum_{j=1}^{k}a_j=\int_1^{k+1}h(x)dx=2\sqrt{k+1}\, \sin(\alpha\log (k+1)).
\]
The function $\sqrt x\sin(\alpha\log x)$ is increasing for $1\le x \le n+1$ because 
$$
(2\sqrt x\sin(\alpha\log x))'=h(x)\geq 0
$$
and consequently
\[
\frac{\sqrt{k+1}}{k^2}\sin(\alpha\log (k+1))\geq \int_k^{k+1}\frac{\sin(\alpha\log x)}{x^{3/2}}dx.
\]
Then
\begin{align*}%\label{eq21}
M_i^*&=2\sum_{k=i}^{n}\frac{\sqrt{k+1}}{k^2}\sin(\alpha\log (k+1))
\geq 2\int_i^{n+1}\frac{\sin(\alpha\log x)}{x^{3/2}}dx\\
&\geq 2\int_{\eta_i}^{n+1}\frac{\sin(\alpha\log x)}{x^{3/2}}dx
=-\frac{4}{1+4\alpha^2}\int_{\eta_i}^{n+1}h'(x)dx=\frac{4}{1+4\alpha^2}h(\eta_i)=\frac{4a_i}{1+4\alpha^2}.
\end{align*}
Thus,
\[
M_i\geq\frac{4}{1+4\alpha^2},\quad i=1,2,...,n,
\]
and 
\[
d_n\geq \min_{1\le i\le n}M_i \geq     \frac{4}{1+4\alpha^2}.
\]
The proof of Theorem \ref{th2} is complete.

%==================================================================
\section{Hardy's  inequality \eqref{eq14}, the zeros of continuous dual Hahn polynomials and the proof of Theorem \ref{Th1.4}}
%==================================================================

In this section we prove Theorem \ref{Th1.4}. 

The continuous dual Hahn polynomials are defined by (see \cite[Section 1.3]{KS1998})
$$
\frac{S_{n}(x^2;a,b,c)}{(a+b)_n (a+c)_n} =  \ _3 F_2 \left( 
\begin{array}{c}
- n,a+ix,a-ix\\
a+b,a+c
\end{array} ; 1 \right)=  \sum_{\nu=0}^{n} \frac{(-n)_{\nu}\,(a+ix)_{\nu}\,(a-ix)_{\nu}}{\nu!\,(a+b)_{\nu}\,(a+c)_{\nu}},
$$
where Pochhammer's symbol is given  by $(\alpha)_\nu=\alpha(\alpha+1)\cdot\ldots\cdot (\alpha+\nu-1)$, $\nu\geq 1$, and $(\alpha)_0:=1$. 
It is clear that each $S_{n}(x^2;a,b,c)$  is a polynomial of degree $2n$ with leading coefficient $(-1)^n$. Moreover, when the parameters $a, b$ and 
$c$ are positive real numbers, $S_{n}(x^2;a,b,c)$ are orthogonal with respect to an absolutely continuous Borel measure (see \cite[(1.3.2) on p. 29]{KS1998}). Hence, the smallest positive zeros $x_{n,1}(a,b,c)$ of $S_{n}(x^2;a,b,c)$ converge to zero 
when $n$ goes to infinity. The following is one of the statements in Theorem 1.1 in \cite{DIGR}:

\begin{THEO}\label{t1.2}
Let $\,d_n\,$ be the smallest possible constant such that 
inequality \eqref{eq14} holds. Then 
\begin{equation}
\label{dnH}
d_n = 4\,  \Bigg(1-\frac{4\, [x_{n,1}(1/2,1/2,1/2)]^2}{1+4\, [x_{n,1}(1/2,1/2,1/2)]^2} \Bigg),
\end{equation}
where $x_{n,1}(1/2,1/2,1/2)$ is the smallest positive zero of $S_{n}(x^2;1/2,1/2,1/2)$.
\end{THEO}

The precise uniform asymptotics for the continuous dual Hahn polynomials was obtained very recently in \cite{Long}. Set
$$
A(z) = \frac{\Gamma(a-z) \Gamma(b-z) \Gamma(c-z)}{\Gamma(1-2z)}.
$$
According to \cite[(80)]{Long},
\begin{equation}
\label{SnAss}
S_n(x^2) = \frac{2 \gamma_{2n+1}}{x \sqrt{w(x)}}\, \cos \left( x \log n + \arg\, A(ix) - \frac{\pi}{2} \right) \left(1 + \mathcal{O}\left( \frac{1}{n} \right) \right),\ \ n\to \infty,
\end{equation}
where $w(x)$ is defined in \cite{Long} and it is always positive. Therefore, we need to calculate the argument of $A(ix)$ for small values of $x$ when $a=b=c=1/2$.  Obviously 
$$
\arg\, A(ix) = \arg\, \frac{[\Gamma(1/2-ix)]^3}{\Gamma(1-2ix)} = 3\, \arg\, \Gamma(1/2-ix) - \arg\, \Gamma(1-2ix)
$$
and
\begin{align*}
  \arg\,\big(\Gamma\big(\tfrac12-ix\big)\big)
  & =  \arcsin \dfrac{ \Im\,\big(\Gamma\big(\tfrac12-ix\big)\big)}
          {\big|\big(\Gamma\big(\tfrac12-ix\big)\big)\big|} = \arcsin \frac{\int_{0}^\infty t^{-1/2} e^{-t} \sin (-x \log t)\, dt}{\big|\big(\Gamma\big(\tfrac12-ix\big)\big)\big|} .
\end{align*}
Since the expansion 
\begin{align*}
\int_{0}^\infty t^{-1/2} e^{-t} \sin (-x \log t)\, dt & = \sum_{k=0}^\infty \frac{(-1)^k\, (-x)^{2k+1}}{(2k+1)!} \int_{0}^\infty t^{-1/2} e^{-t} (\log t)^{2k+1}\, dt\\
& = \sum_{k=0}^\infty \frac{(-1)^{k+1}\, \Gamma^{(2k+1)}(1/2)}{(2k+1)!}\  x^{2k+1}
\end{align*}
certainly holds for $x\in(-1/2,1/2)$ because the radius of convergence of the series is exactly $1/2$, then 
$$
\frac{\Im\,\big(\Gamma\big(\tfrac12-ix\big)\big)}{\big|\big(\Gamma\big(\tfrac12-ix\big)\big)\big|}  \sim \frac{-\, \Gamma^{\prime}(1/2)}{\Gamma(1/2)}\,   x =  (\gamma +\log 4)\ x
$$
for all sufficiently small $x$. Thus, 
$$
\arg\,\big(\Gamma\big(\tfrac12-ix\big)\big) \sim \arcsin \big( (\gamma +\log 4)\ x \big) \sim (\gamma +\log 4)\ x\ \ \ \mathrm{as}\ \ x \to 0.
$$

Similar reasonings show that 
\begin{align*}
  \arg\,\big(\Gamma\big(1-2ix\big)\big)
  & =  \arcsin \dfrac{ \Im\,\big(\Gamma\big(1-2ix\big)\big)}
          {\big|\big(\Gamma\big(1-2ix\big)\big)\big|} = \arcsin \frac{\int_{0}^\infty e^{-t} \sin (-2x \log t)\, dt}{\big|\big(\Gamma\big(1-2ix\big)\big)\big|},
\end{align*}

\begin{align*}
\int_{0}^\infty e^{-t} \sin (-2 x \log t)\, dt & = \sum_{k=0}^\infty \frac{(-1)^{k+1}\, \Gamma^{(2k+1)}(1)}{(2k+1)!}\  (2x)^{2k+1}
\end{align*}

and
$$
\arg\,\big(\Gamma\big(1-2ix\big)\big) \sim \frac{-\, \Gamma^{\prime}(1)}{\Gamma(1)}\,   (2x) =  2\, \gamma\, x\ \ \ \mathrm{as}\ \ x \to 0.
$$
Hence 
$$
\arg\, A(ix) \sim 3\, (\gamma +\log 4)\, x - 2\, \gamma\, x = (\gamma + \log 64)\, x,\ \ \ x \to 0.
$$
Since we are interested in the smallest zero, then indeed $x_{n,1}(1/2,1/2,1/2)$ converges to zero when $n$ goes to infinity. Therefore, we use the latter approximation 
of  $arg\, A(ix)$. Thus, for the argument $u(x)$ of the cosine in (\ref{SnAss}) we obtain
$$
u(x)= x \log n + (\gamma + \log 64)\, x -\frac{\pi}{2},
$$
so that $x_{n,1}(1/2,1/2,1/2)$ is asymptotically equal to the smallest positive zero of $\cos u(x)=0$, that is, to the solution of $u(x)=\pi/2$. In other words,
$$
 x_{n,1}(1/2,1/2,1/2) \sim \frac{\pi}{\gamma + \log 64 + \log n}, \ \ \mathrm{as}\ \ n \to \infty.
$$ 
The latter, together with (\ref{dnH}), yields
$$
d_n \sim 4 - \frac{16 \, \pi^2}{4\, \pi^2 + (\gamma + \log 64 + \log n)^2},\ \ \mathrm{as}\ \ n \to \infty,
$$
which is exactly (\ref{asdn}). 

The explicit form of the coefficients $c_k$ in (\ref{ck}) follows immediately after straightforward manipulations with the error term on the right-hand side of the latter expression and the generating function of   
the Chebyshev polynomials of the second kind 
$$
\frac{1}{1-2\, t\, y + y^2} = \sum_{j=0}^\infty\, U_j(t)\, y ^j.
$$
Indeed, setting 
$$
y_n = \frac{\left(4 \pi^2 + (\gamma+6\log2)^2\right)^{1/2}}{\log n}\ \ \ \mathrm{and}\ \ \ \tilde{t}= - \frac{\gamma+6\log2}{\left(4 \pi^2 + (\gamma+6\log2)^2\right)^{1/2}},
$$
the error term of (\ref{asdn}) becomes
$$
 \frac{16 \, \pi^2}{4\, \pi^2 + (\gamma + \log 64 + \log n)^2} = \frac{16 \, \pi^2}{[\log n]^2}  \frac{1}{1-2\, \tilde{t}\, y_n + y_n^2} = \frac{16 \, \pi^2}{[\log n]^2}\  \sum_{j=0}^\infty\, U_j(\tilde{t})\, y_n ^j. 
 $$
 Since $|\tilde{t}| <1$ and all Chebyshev polynomials of the second kind obey the inequality $|U_j(t)|\leq (1-t^2)^{-1/2}$ for $t\in (-1,1)$, the latter series is absolutely convergent for $|y_n| <1$. 
 The same argument shows that the error term in (\ref{asexp}) is indeed $O((\log n)^{-m-1})$. 
 
 \section*{Acknowledgments}
We thank F. Stampach for calling our attention to Wilf's and his own contributions.

\end{document}